\documentclass{article}

\usepackage{graphicx}
\usepackage{amsmath, amsthm}
\usepackage{amsfonts}
\usepackage{amssymb}
\usepackage{url}
\usepackage{hyperref}  
\usepackage{verbatim}
\usepackage{color}
\definecolor{red}{rgb}{1,0,0}

\definecolor{blu}{rgb}{0,0,1}

\usepackage[mathlines]{lineno}
\usepackage{dsfont} 

\def\noi{\noindent}
\def\mtx#1{\begin{bmatrix} #1 \end{bmatrix}}

\newcommand{\nul}{\operatorname{null}}

\newcommand{\R}{\mathbb{R}}

\newcommand{\C}{\mathbb{C}}

\newcommand{\Rnn}{\R^{n\times n}}

\newcommand{\bx}{{\bf x}}
\newcommand{\by}{{\bf y}}

\newcommand{\bone}{\ensuremath{\mathds{1}}} 

\newcommand{\pd}{\gamma_P}

\newcommand{\Z}{\operatorname{Z}}
\newcommand{\mult}{\operatorname{mult}}
\newcommand{\clf}{\mathcal{F}}
\newcommand{\E}{\mathbb{E}}

\newcommand{\G}{\mathcal{G}}
\newcommand{\SG}{\mathcal{S}(G)}

\newcommand{\sym}{\mathcal{S}}
\newcommand{\mat}{\mathcal{M}_0}
\newcommand{\mr}{\operatorname{mr}}

\newcommand{\M}{\operatorname{M}}

\newcommand{\Zm}{\operatorname{Z}^-}

\newcommand{\rank}{\operatorname{rank}}

\newtheorem{thm}{Theorem}[section]
\newtheorem{cor}[thm]{Corollary}

\newtheorem{prop}[thm]{Proposition}

\newtheorem{obs}[thm]{Observation}

\theoremstyle{definition}
\newtheorem{rem}[thm]{Remark}

\theoremstyle{definition}

\theoremstyle{definition}
\newtheorem{ex}[thm]{Example}

\newcommand{\x}{\times}

\newcommand{\bit}{\begin{itemize}}
\newcommand{\eit}{\end{itemize}}
\newcommand{\ben}{\begin{enumerate}}
\newcommand{\een}{\end{enumerate}}
\newcommand{\beq}{\begin{equation}}
\newcommand{\eeq}{\end{equation}}
\newcommand{\bea}{\begin{eqnarray*}}
\newcommand{\eea}{\end{eqnarray*}}
\newcommand{\bpf}{\begin{proof}}
\newcommand{\epf}{\end{proof}\ms}
\newcommand{\ms}{\medskip}

\newcommand{\cp}{\, \Box\,}
\newcommand{\lc}{\left\lceil}
\newcommand{\rc}{\right\rceil}
\newcommand{\lf}{\left\lfloor}
\newcommand{\rf}{\right\rfloor}

\title{Zero forcing and power domination for graph products}

\author{Katherine F. Benson\thanks{Department of Mathematics and Physics, Westminster College, Fulton, MO 65251, USA (katie.benson@westminster-mo.edu).} \and Daniela Ferrero\thanks{Department of Mathematics, Texas State University, San Marcos, TX 78666, USA (dferrero@txstate.edu).} \and Mary Flagg\thanks{Department of Mathematics, Computer Science and Cooperative Engineering, University of St. Thomas, 3800 Montrose, Houston, TX 77006, USA (flaggm@stthom.edu).} \and Veronika Furst\thanks{Department of Mathematics, Fort Lewis College, Durango, CO 81301,
USA (furst$\_$v@fortlewis.edu).} \and Leslie Hogben\thanks{Department of Mathematics, Iowa State University,
Ames, IA 50011, USA (hogben@iastate.edu) and American Institute of Mathematics, 600 E. Brokaw Road, San Jose, CA 95112, USA
(hogben@aimath.org).} \and Violeta Vasilevska\thanks{Department of Mathematics, Utah Valley University, Orem, UT, 84058, USA (Violeta.Vasilevska@uvu.edu).} \and Brian Wissman\thanks{Natural Sciences Division, University of Hawaii at Hilo, Hilo, HI 96720-4091, USA (wissman@hawaii.edu).}}

\begin{document}

\maketitle


\vspace{-10pt}
\begin{abstract} 

The power domination number arises from the monitoring of electrical networks and its determination is an important problem.  Upper bounds for power domination numbers can be obtained by constructions.  Lower bounds for the power domination number of several families of graphs are known, but they usually arise from specific properties of each family and the methods do not generalize.  In this paper we  exploit the relationship between power domination and zero forcing to obtain the first general lower bound for the power domination number.  
We apply this bound to obtain results for both the power domination of tensor products and the zero-forcing number of lexicographic products of graphs.  We also establish results for the zero forcing number of  tensor products and Cartesian products  of graphs.
\end{abstract}

\noi {\bf Keywords} power domination, zero forcing, maximum nullity, minimum rank, tensor product, lexicographic product, Cartesian product 

\noi{\bf AMS subject classification} 05C69, 05C50 


\section{Introduction}\label{sintro}

 Electric power companies need to monitor the state of their networks
continuously. 
One method of monitoring a network is to place Phase Measurement Units (PMUs) at selected
locations in the system, called electrical nodes or buses, where transmission lines, loads, and generators are connected.  A PMU placed at an electrical node measures the voltage at the node and all current phasors at the node \cite{BMBA93}; it also provides these measurements at other vertices or edges according to certain propagation rules (see Section \ref{ssdefnot} for the formal definitions of this and other terms).  The placement of PMUs at all nodes of a network is a trivial solution to the monitoring problem. Because of the cost of a PMU, the trivial solution is not feasible, and it is
important to minimize the number of PMUs used while maintaining the
ability to observe the entire system.

This problem was first studied in terms of graphs by Haynes et al.~in
\cite{HHHH02}. Indeed, an electric power network can be modeled by a
graph where the vertices represent the electric nodes and the edges
are associated with the transmission lines joining two electrical
nodes. In this model, the power domination problem in graphs
consists of finding a minimum set of vertices from where the entire
graph can be observed according to certain rules. In terms of the
physical network, such a minimum set of  vertices will provide the locations where
the PMUs should be placed in order to monitor the entire graph at
 minimum cost. Since its introduction in \cite{HHHH02}, the power domination number and variations have generated considerable interest (see, for example, \cite{BH, CDMR12, Deanetal11, DMKS, DH06, L14, SRRGW15}).

As was pointed out in \cite{Deanetal11}, a careful examination of the definition of power domination leads naturally to the study of zero forcing.  The zero forcing number was introduced in \cite{AIM08}  as an upper bound for the maximum nullity of real symmetric matrices whose nonzero pattern of off-diagonal entries is described by a given graph, and independently by mathematical physicists studying control of quantum systems.

In Section \ref{spdzf} we use the connection between the two definitions established in \cite{Deanetal11} to obtain the only known general lower bound for the power domination number. Note that in \cite{L14}  the author claimed to have obtained the first general lower bound for the power domination number, but a family of counterexamples to his claim was given in \cite{FHKY15}.  Then we use our lower bound to prove results for both power domination and zero forcing.  Some of these applications require additional knowledge of zero forcing numbers of tensor products, which we establish in Section \ref{Ztensor}, via a new upper bound on the zero forcing number of the tensor product of a complete graph with another graph.  
The remainder of this introduction contains formal definitions of power domination and zero forcing, graph terminology, and matrix terminology.


\subsection{Power domination and zero forcing definitions}\label{ssdefnot}

A  {\em graph} $G=(V,E)$  is an ordered pair formed by a finite nonempty set of {\em  vertices}
$V=V(G)$ and a set of {\em edges} $E=E(G)$ containing unordered pairs
of distinct vertices (that is, all graphs are simple and undirected).  The {\em order} of $G$ is denoted by $|G|:=|V(G)|$. We say the vertices $u$ and $v$ are {\em adjacent} or are {\em neighbors}, and write $u\sim v$, if $\{u,v\} \in E$.   For any vertex $v \in V$, the {\em neighborhood} of $v$
is the set $N(v) = \{u \in V: u\sim v\}$  (or $N_G(v)$ if $G$ is not clear from context), and the
{\em closed neighborhood} of $v$ is the set $N[v] = N(v) \cup \{ v \}$.
Similarly, for any set of vertices $S$, $N(S) = \cup_{v \in S} N(v)$
and $N[S] = \cup_{v \in  S} N[v]$.

A vertex $v$ in a graph $G$ is said to \emph{dominate} itself and all of its neighbors in $G$.  A set of vertices $S$ is a \emph{dominating set} of $G$ if every vertex of $G$ is dominated by a vertex in $S$.  The minimum cardinality of a dominating set is the \emph{domination number} of $G$ and is denoted by $\gamma (G)$.

In \cite{HHHH02}  the authors introduced the related concept of power domination by presenting propagation rules in terms of vertices and edges in a graph. In this paper we will use a simplified version of the propagation rules that is equivalent to the original, as shown in   \cite{BH}. For a  set $S$ of vertices in a graph $G$, define  $PD(S)\subseteq V(G)$ recursively:

\ben
\item $PD(S):=N[S]= S\cup N(S)$.
\item While there exists $ v\in PD(S)$ such that $|N(v)\cap (V(G)\setminus PD(S))|=1$: 

$PD(S):=PD(S)\cup N(v)$.
\een

We say that a set $S\subseteq V(G)$ is a {\em power dominating set} of  a graph $G$ if at the end of the process above $PD(S)=V(G)$. A {\em minimum power dominating set} is a power dominating set of minimum cardinality, and the {\em power domination number}  $\pd (G)$ of  $G$ is the cardinality of a minimum power dominating set.

The concept of zero forcing can be explained via a coloring game on the vertices of $G$.  The {\em color change rule} is: If $u$ is a blue vertex and exactly one  neighbor $w$ of $u$ is white, then change the color of $w$ to blue.  We say $u$ forces $w$ and denote this  by $u\to w$.
A {\em zero forcing set} for  $G$ is a subset of vertices $B$ such that when the vertices in $B$ are colored blue and the remaining vertices are colored white initially,  repeated application of the color change rule can color all vertices of $G$ blue.  A {\em minimum zero forcing set} is a zero forcing set of minimum cardinality, and the {\em zero forcing number} $\Z(G)$ of  $G$ is the cardinality of a minimum zero forcing set.
  The next observation is the key relationship between the two concepts.

 \begin{obs}\label{keyZpd} {\rm \cite{Deanetal11}}  The power domination process on a graph $G$ can be described as choosing a set $S\subseteq V(G)$ and applying the zero forcing process to the closed neighborhood $N[S]$ of $S$.  The set $S$ is a power dominating set of $G$ if and only if $N[S]$ is a zero forcing set for $G$.  
 \end{obs}

The {\em degree} of a vertex $v$, denoted by $\deg v$, is the
cardinality of the set $N(v)$. The {\em maximum} and {\em minimum degree} of $G$ are defined as $\Delta(G) = \max \{\deg v:v\in V\}$ and $\delta(G) = \min \{\deg v:v\in V\}$, respectively. A graph $G$ is {\em regular} if $\delta(G)=\Delta(G)$. 

The next observation is well known (and immediate since the color change rule cannot be applied in $G$ without at least $\delta(G)$ blue vertices).

 \begin{obs}\label{deltaZ} For every graph $G$, $\delta(G)\le\Z(G)$.
 \end{obs}


\subsection{Graph definitions and notation}\label{ssdefnot}

Let $n$ be a positive integer. The {\em path} of order $n$ is the graph $P_n$
with   $V(P_n) = \{ x_i: 1 \leq i \leq n \}$ and   $E(P_n) = \{ \{  x_i,  x_{i+1} \}:1 \leq i \leq n-1 \}$.  If $n\geq 3$, the {\em cycle} of order $n$ is the graph $C_n$ with  $V(C_n) = \{ x_i:
1 \leq i \leq n \}$ and  $E(C_n) = \left\{ \left\{  x_i,  x_{i+1} \right\}: 1 \leq i \leq n-1 \right\}\cup\left\{\{x_n,x_1\}\right\}$.  The {\em complete graph} of order $n$ is the graph $K_n$
with  $V(K_n) = \{ x_i: 1 \leq i \leq n \}$ and $E(K_n) = \{ \{  x_i,  x_j \}: 1 \leq i<j \leq n \}$.
%

Let $G = (V(G), E(G))$ and $H = (V(H), E(H))$ be disjoint graphs. All of the following products of $G$ and $H$ have vertex set $V(G)\x V(H)$.  The {\em tensor product}  (also called the {\em direct product}) of $G$ and $H$ is denoted by $G\times H$; a vertex $(g,h)$ is adjacent to  a vertex $(g',h')$ in $G\times H$ if  $\{g,g'\} \in E(G)$ and $\{h,h'\} \in E(H)$.  The \emph{Cartesian product} of $G$ and $H$ is denoted by $G \cp H$; two  vertices $(g,h)$ and $(g',h')$ are adjacent in $G \cp H$ if either (1) $g=g'$ and $\{h,h'\} \in E(H)$, or (2) $h=h'$ and $\{g,g'\} \in E(G)$. 
The {\em lexicographic product} of $G$ and $H$ is denoted by  $G*H$; two vertices $(g,h)$ and $(g', h')$ are adjacent in $G*H$ if either (1) $\{g,g'\} \in E(G)$, or (2) $g=g'$ and $\{h,h'\} \in E(H)$.  Note that $H\x G\cong G\x H$ and $H\cp G\cong G\cp H$, whereas $H * G$ need not be isomorphic to $G*H$.

For a graph $G$ with no edges, $\Z(G)=\pd(G)=\gamma(G)=|G|$, so we focus our attention on graphs with edges. 
In the case of the tensor product $G\times H$, this means we assume $|G|, |H| \geq 2$.


\subsection{Matrix definitions and notation}\label{ssdefnot}

 Let $S_n(\R)$ denote the set of all $n\times n$ real symmetric matrices. For $A = [a_{ij}] \in S_n(\R)$,
the {\em graph} of $A$,  denoted by $\G(A)$, is the graph with vertices $\{1, \ldots , n\}$ and edges
$\{\{i, j\}: a_{ij}\not= 0,\,  1 \leq i < j \leq n\}$. More generally, the graph of $A$ is defined for any matrix that is {\em combinatorially symmetric}, i.e., $a_{ij}=0$ if and only if $a_{ji}=0$.  Note that the diagonal of $A$ is ignored in determining
$\G(A)$. The set of symmetric matrices described by a graph $G$ of order $n$ is defined as  $\sym(G)=\{A \in S_n(\R) : \G(A) = G\}$.
The {\em maximum nullity} of $G$ is $\M(G) = \max \{\nul A: A \in \sym(G)\}$, and the {\em minimum rank} of $G$ is $\mr(G) = \min \{{\rm rank}\, A : A \in S(G)\}$; clearly $\M(G)+\mr(G)=|G|$. The term `zero forcing' comes from using the forcing process to force zeros in a null vector of a matrix $A\in\sym(G)$, implying  $\M(G)\le\Z(G)$ \cite{AIM08}.

A standard way to construct matrices of maximum nullity for  a Cartesian product or a  tensor product of graphs is to use the {\em Kronecker} or {\em tensor} product of matrices.  Let $A$ be an $n \times n$ real matrix and $B$ be an $m \times m$ real matrix.  Then $A \otimes B$ is the $n \times n$ block matrix whose $ij$th block is the $m\times m$ matrix  $a_{ij}B$. It is known that $(A\otimes B)^T=A^T\otimes B^T$ and $\rank( A\otimes B)=(\rank A)(\rank B)$.  
 If  $A\in \sym(G)$, $B\in\sym(H)$, $|G|=n$, and $|H|=m$, then $A \otimes I_m - I_n \otimes B \in \sym(G\cp H)$. 
If $\bx$ is an eigenvector of 
$A$ for eigenvalue $\lambda$ and $\by$ is an eigenvector of 
$B$ for eigenvalue $\mu$, then  $\bx\otimes\by$ is an eigenvector of $A \otimes I_m - I_n \otimes B$ for eigenvalue $\lambda-\mu$.  Since a real symmetric matrix has a full set of eigenvectors, the multiplicity of $\lambda-\mu$ is   $\mult_A(\lambda)\, \mult_B(\mu)$    \cite[Observation 3.5]{AIM08}.
If $A\in \sym(G)$ and $B\in\sym(H)$ and the diagonal entries of $A$ and $B$ are all zero,  then  $A\otimes B\in\sym(G\x H)$. Define $\mat(G)=\{A \in \Rnn : \G(A) = G \mbox{ and $a_{ii}=0$ for } i=1,\dots,n\}$; in contrast to a matrix in $\sym(G)$, a matrix in $\mat(G)$ need not be symmetric but must have a zero diagonal and be {combinatorially symmetric}.  If $A\in\mat(G)$ and $B\in\mat(H)$, then  $A\otimes B\in\mat(G\x H)$. 


\section{Zero forcing for graph products}\label{szfs}

In this section we develop a tool for the zero forcing number of tensor products of graphs and apply it to compute the tensor products of complete graphs with paths and with cycles.  We also compute the zero forcing number and maximum nullity of the Cartesian product of two cycles.  


\subsection{Tensor products}

For tensor products of graphs we use not only  standard zero forcing, but also skew zero forcing \cite{IMA} defined by the \emph{skew color change rule}: If a vertex $v$ of $G$ has exactly one white neighbor $w$, then $v$ forces $w$ to change color  to blue;  if $v$ is white  when it forces, the force $v\to w$ is called  a {\em white vertex force}.   The \emph{skew zero forcing number}, $\Zm(G)$, is the minimum cardinality of a {\em skew zero forcing set}, i.e., a (possibly empty) set of blue vertices that can color all vertices blue using the skew color change rule. Note that we are using skew zero forcing as a tool, and it is not directly connected to power domination.  
For either a   standard or skew zero forcing set $B$, color all the vertices of $B$ blue and then perform zero forcing, listing the forces in the order in which they were performed.  This list is a {\em chronological list of forces of $B$} and is denoted by $\clf$.

For each vertex $g \in V(G)$, define the set $U_{g}=\{(g,h): h\in V(H)\}$ in $G\x H$.  We say vertices $(g,h)$ and $(g',h)$ are {\em associates} in $G\x H$ if $g\sim g'$ in $G$; associates are not adjacent in $G\x H$.

\begin{thm}\label{Ztensor}
Let $G$ be a graph  and $n\geq 4$.  Then \[\Z(G\times K_n)\le (n-2)|G| + 2\Zm(G).\]
\end{thm}

\bpf    For $g\sim g'$ in $G$,   $(g,i)\sim (g',j)$ in $G\x K_n$ for all $j\ne i\in V(K_n)$.  Choose a minimum skew zero forcing set $B$ for $G$ and a chronological list of forces $\clf$ of $B$ and denote the $k$th force by $g_k\to w_k$ (many vertices receive two labels, e.g., $g_k=w_\ell$).  We describe how to choose $ (n-2)|G| + 2|B|$  vertices  to obtain a zero forcing set $\hat B$ for $G\times K_n$.  For $g\in B$, let $ \hat B \supset U_g$.  For $v\not\in B$, place $n-2$ vertices of $U_v$ in $\hat B$; the selection of these vertices is determined when $v$ is forced  in $G$ or when $v$ performs a force in $G$, whichever comes first. 

Consider the $k$th force in $\clf$ ($k=1$ is permitted).  
  Suppose $g_k\to w_k$ is not a white vertex force.  If no vertices of  $U_{w_k}$ are in $\hat B$ yet, then arbitrarily choose $n-2$ vertices in $U_{w_k}$ to place in $\hat B$; otherwise, no additional vertices  are placed in $\hat B$. 
Now suppose $g_k\to w_k$ is a white vertex force.  Clearly, $w_k\not\in B$ and $w_k$ has not been forced previously in $G$. The only force ${w_k}$ could have performed in $G$ would be  $w_k\to {g_k}$, in which case $g_k\to w_k$ would  not be a white vertex force.  Thus, no vertices in $U_{w_k}$ have been previously placed in $\hat B$.  Also, no vertices of $U_{g_k}$ have been previously placed in $\hat B$.  Since $n\geq 4$, it is possible to choose $n-2$ vertices in each of $U_{g_k}$ and $U_{w_k}$, and place in $\hat B$, in such a way that for every pair of associated vertices in $U_{g_k}$ and $U_{w_k}$, at least one member of the pair is placed in $\hat B$.  By construction, $|\hat B| = |B|n + (|G|-|B|)(n-2)$.

The zero forcing process in $G\times K_n$ now follows the chronological list of forces $\clf$. 
At stage $k$ (before performing the $k$th force), we assume every vertex in $U_v$ is blue for every $v$ such that $v\in B$ or $v=w_\ell$ with $\ell <k$. Given the force $g_k\to w_k$ in $\clf$, the two (blue) vertices in $U_{g_k}$ that are associated with the two white vertices in $U_{w_k}$ can each force the other's associate in $U_{w_k}$, so $U_{w_k}$ is now entirely blue. Thus, all sets $U_{w_k}$ will be turned entirely blue. \epf

\begin{figure}[!ht]
\begin{center}
\scalebox{.4}{\includegraphics{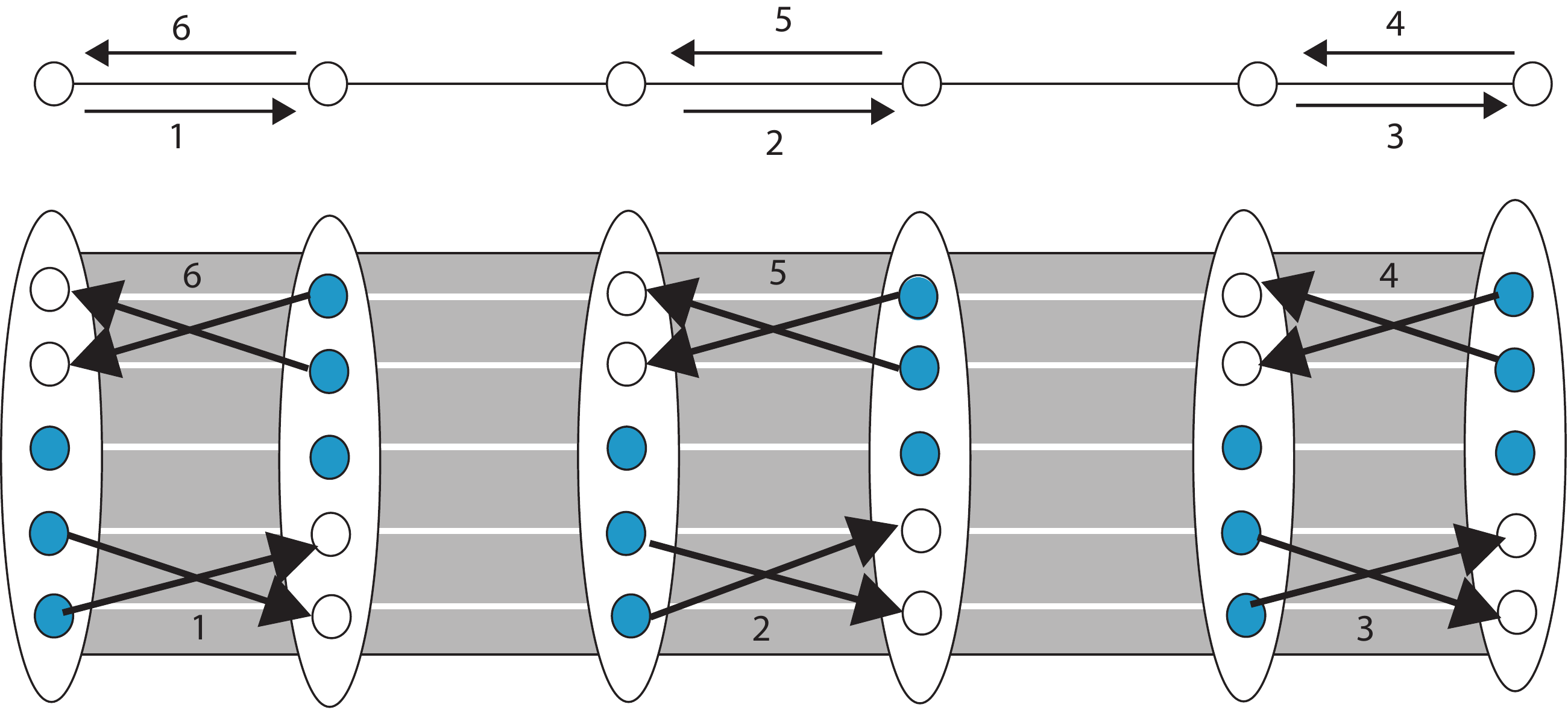}} \vspace{-8pt}
\caption{The skew zero forcing process on $P_6$ and the analogous zero forcing process on $P_{6}\x K_5$ are illustrated.  In the schematic diagram of $P_{6}\x K_5$, the gray areas indicate that all possible edges are present except for the non-edges marked as white lines.  \vspace{-15pt}} \label{f:ZPellxKn}
 \end{center}
 \end{figure}

The forcing process used in the proof of Theorem \ref{Ztensor} is illustrated   in Figure \ref{f:ZPellxKn}, where skew forcing is shown on $P_6$ and standard zero forcing is shown on $P_6\x K_5$; the number of each force in a chronological list of forces of $P_6$ is also shown. The first half of the forces in $P_6$ (forces 1, 2, and 3) are white vertex forces and the second half are not.

Note that the bound in Theorem \ref{Ztensor} need not be valid for $n=3$, as shown in the next example.

\begin{ex}\label{3sun-counterex}  Let $H_3$ denote the 3-sun shown in Figure \ref{f:3sun-counterex} and consider $H_3\x K_3$. Suppose $\hat B$ is a zero forcing set for $H_3\x K_3$ of cardinality six.   Observe that $\hat B$ must contain at least one vertex in  $U_g$ for each $g\in V(H_3)$,  or no vertices in $U_g$ could ever be colored blue, so necessarily $\hat B$ contains exactly one vertex of each  $U_g$.  In $H_3\x K_3$ there are three sets $U_g$ that contain vertices of degree 2 (corresponding to the three vertices of degree 1 in $H_3$) and three sets of vertices of degree 6.  With appropriately staggered choices of which vertex of $U_g$ is  in $\hat B$, each of the three blue vertices of degree 2 can force one degree 6 vertex blue.  Now the only blue vertices that have not yet forced all have degree 6.  Each has (at least) one white neighbor of degree 2.  In order for such vertex to perform a force, it must have no other white neighbors.  That is, all four of its degree 6 neighbors must be blue.  This implies  two degree 6 white vertices must be associates, preventing any further forcing. Thus $\Z(H_3\x K_3)\ge 7>6 = (3-2)6+0$ and $\Zm(H_3)=0$.
(The zero forcing number of this example was originally found by use of the {\em Sage} zero forcing software \cite{sagedata}, which found a zero forcing set of order 7 and determined no smaller ones exist.)

\begin{figure}[!ht]
\begin{center}
\scalebox{.4}{\includegraphics{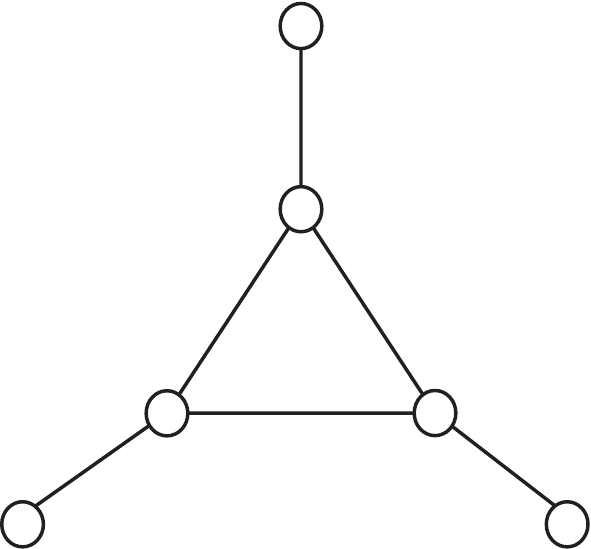}} \vspace{-8pt}
\caption{The 3-sun  \vspace{-25pt}} \label{f:3sun-counterex}
 \end{center}
 \end{figure}

\end{ex}

\vspace{-15pt} We apply Theorem \ref{Ztensor}  to the tensor product of a path and a complete graph.  The case of odd paths has already been  done:

\begin{thm}{\rm \cite[Theorem 15]{HCY10}} \label{HCYthm}
If  $t\geq 1$ is odd and $n\geq 2$, then $M(P_t\times K_n)=Z(P_t\times K_n)=(n-2)t+2$.
\end{thm}

The method used to prove  \cite[Theorem 15]{HCY10} is the standard one of exhibiting a matrix  and a zero forcing set with   cardinality equal to the nullity of the matrix.  For  even paths, we  use the same matrix to establish a lower bound on maximum nullity.  For $n\ge 4$, Theorem \ref{Ztensor} gives an equal upper bound on the zero forcing number;   a specific zero forcing set is exhibited for $n=3$.

\begin{thm}\label{evenPtKn}
If  $t\geq 2$ is even and $n\geq 3$, then \[\M(P_{t}\times K_n)=\Z(P_{t}\times K_n)= (n-2)t.\]
\end{thm}

\bpf  Define the $n$-vectors  $\bone=[1,1,\dots,1]^T$ and $\by=[1,2,\dots,n]^T$, and define $A_n=[\bone\, \by]\mtx{0 & 1\\-1 & 0}\mtx{\bone^T\\\by^T}$.  Then $\rank A_n=2$, $A_n^T=-A_n$, and $A_n\in\mat(K_n)$.  Define $B_t$ to be a $t\x t$ skew adjacency matrix (a tridiagonal matrix with 1s on the first superdiagonal, 0s on the main diagonal, and $-1$s on the first subdiagonal); note that since $t$ is even, $\rank B_{t}=t$.  Then $B_{t}\otimes A_n\in\sym(P_t\x K_n)$ and $\rank(B_{t}\otimes A_n)=2t$.  Thus $(n-2)t=t n-2t\le \M(P_{t}\times K_n) \le \Z(P_{t}\times K_n)$.  
For $n\ge 4$, $\Z(P_{t}\times K_n)\le (n-2)t$ 
by Theorem \ref{Ztensor}, since $\Zm(P_t)=0$ for $t$ even \cite{IMA}.

Now assume $n=3$.  The forcing order in $P_t \times K_3$ is slightly more complicated than the one in Theorem \ref{Ztensor}, since initially only a single vertex can force at a time.  Label the vertices  of  $P_{t}\times K_3$ as ordered pairs $(r,s)$ with $1\le r\le t$ and $1\le s\le 3$.  Define $B= \{(2i-1,3),(2i,1):i=1,\dots, \frac t 2\}$.
 First, $(1,3)$ forces $(2,2)$.  
 Continue in increasing order of sets, so $(2i-1,3)$ forces $(2i,2)$ for  $i=1,\dots, \frac t 2$.  Then the process is repeated in reverse order, starting at $2i=t$.  Now, $(t,1)$ and $(t,2)$ are both blue with one white neighbor each, so $(t,1)$ forces $(t-1,2)$ and $(t,2)$ forces $(t-1,1)$.  Continue in decreasing order, so $(2i,1)$ forces $(2i-1, 2)$ and $(2i,2)$ forces $(2i-1,1)$ for  $i= \frac t 2$ down to $i=1$, turning all odd-numbered sets all blue.  Finally, in increasing order again, $(2i-1,2)$ forces $(2i, 3)$ for  $i=1,\dots, \frac t 2$.
\epf

Observe that the formula in Theorem \ref{evenPtKn} fails for $n = 2$:  $P_t\x K_2$ is the disjoint union of two copies of $P_2$, so $\M(P_t\x K_2)=\Z(P_t\x K_2)=2\ne(2-2)t$.  

\begin{thm}\label{CtKn}
If  $n, t\geq 3$, then \[\M(C_{t}\times K_n)=\Z(C_{t}\times K_n)=\left\{\begin{array}{ll}
(n-2)t+2 &  \mbox{ if $t$ is odd,}\\
(n-2)t+4 &  \mbox{ if $t$ is even.}\end{array}\right. \]
\end{thm}

\bpf  Let the matrix $A_n$ be as defined in the proof of Theorem \ref{evenPtKn}, so $\rank A_n=2$, $A_n^T=-A_n$, and $A_n\in\mat(K_n)$.  Define $B_t$ to be the $t\x t$ skew-adjacency matrix (with 1s in one cyclic direction and $-1$s in the other).  It is easy to verify that $\det B_t=0$,   $\rank B_{t}=t-1$ for $t$ odd, and $\rank B_{t}=t-2$ for $t$ even. 
Then $B_{t}\otimes A_n\in\sym(C_t\x K_n)$ and 
\[\rank(B_{t}\otimes A_n)=\left\{\begin{array}{ll}
2t-2 &  \mbox{ if $t$ is odd,}\\
2t-4 &  \mbox{ if $t$ is even.}\end{array}\right. \]

Since $ \M(C_{t}\times K_n) \le \Z(C_{t}\times K_n)$, it suffices to exhibit a zero forcing set of cardinality $(n-2)t+2$ for $t$ odd and $(n-2)t+4$ for $t$ even.  Color all the vertices in $U_{t}$ blue for $t$ odd, or all the vertices in $U_{t-1}$ and in $U_{t}$ blue for $t$ even. 
We  now consider the graph obtained by deleting these blue vertices,  i.e., $P_{t-1}\x K_n$  or $P_{t-2}\x K_n$,  respectively, both having the form of a tensor product of an even path with a complete graph.  We construct a minimum zero forcing set $\hat B'$ (of cardinality $(n-2)(t-1)$ or $(n-2)(t-2)$, respectively)  and  perform zero forcing  for the tensor product of an even path and a complete graph as in  Theorem \ref{evenPtKn}. The zero forcing set for $C_{t}\times K_n$ is $\hat B' \cup U_t$ or $\hat B' \cup U_t\cup U_{t-1}$, respectively. 
\epf


\subsection{Cartesian products}

Next we determine the zero forcing number and maximum nullity of the Cartesian product of two cycles.  

\begin{thm}\label{CnCm} For $m\ge n\ge 3$, 
\[\M(C_n\cp C_m)=\Z(C_n\cp C_m)= \left\{ \begin{array}{ll}
 2n-1& \mbox { if $m=n$ and $n$ is odd,}\\[1mm]
 2n & \mbox{ otherwise.}
 \end{array} \right.\]
\end{thm}
\bpf For $m= n\ge 3$, by \cite[Theorem 2.18]{UOM} $\M(C_n\cp C_n)=\Z(C_n\cp C_n)=n+2\lf \frac n 2\rf$, so $\M(C_n\cp C_n)=\Z(C_n\cp C_n)=2n-1$ for $n$ odd and  $\M(C_n\cp C_n)=\Z(C_n\cp C_n)=2n$ for $n$ even. 

So assume $m>n\ge 3$.  
 It is easy to see that the vertices of two consecutive cycles $C_n$ form a zero forcing set, so $\Z(C_n\cp C_m)\le 2n$.  To complete the proof we construct a matrix in $\sym(C_n\cp C_m)$ with nullity $2n$, so $2n\le \M(C_n\cp C_m)\le \Z(C_n\cp C_m)\le 2n$.

 Let  
 $k=\lc \frac n 2 \rc$.  Let $A$ be the matrix obtained from the adjacency matrix of $C_n$ by changing  one pair of  symmetrically placed entries from 1 to $-1$.  Then the   distinct eigenvalues of $A$ are $\mu_i=2\cos \frac{\pi (2i-1)}n,$ $i= 1, \dots , {k}$, each with multiplicity 2 except  $\mu_k=-2$, which has  multiplicity 1 when  $n$ is odd \cite[Theorem 3.8]{AIM08}.      Assuming that there exists a matrix $B\in\sym(C_m)$ such that $\mu_i$ is an eigenvalue of $B$ with  multiplicity 2  for $i=1,\dots,k$, it follows that $A \otimes I_m - I_n \otimes B$ has eigenvalue zero with multiplicity $2n$,  because every eigenvalue of $A$ has a corresponding eigenvalue of $B$ with multiplicity 2.

It remains to establish the existence of a matrix $B\in\sym(C_m)$ such that $\mu_i$ is an eigenvalue of  $B$ with  multiplicity 2  for $i=1,\dots,k$.
In Ferguson \cite[Theorem 4.3]{F80} it is shown that for   any set of $m$ real numbers $\lambda_1>\lambda_2\ge \lambda_3>\lambda_4\ge\lambda_5>\dots$, there is a matrix $B\in\sym(C_m)$ having eigenvalues $\lambda_1,\dots,\lambda_m$.  Thus for $m$ odd or $m\ge n+2$ we can choose $\lambda_{2i}= \lambda_{2i+1}=\mu_i$. 
 Hall  \cite{newpaper} showed that  for $m=2k$ and any $\lambda_1>\lambda_2>\dots>\lambda_k$ there is a  matrix $B\in\sym(C_m)$ with $\mult_B(\lambda_i)=2$ for $i=1,\dots,k$.  
  \epf


\section{Zero forcing lower bound for power domination number}\label{spdzf}

 The power domination number of several families of graphs has been determined using a two-step process: finding an upper bound and a lower bound. The upper bound is usually obtained by providing a pattern to construct a set, together with a proof that the constructed set is a power dominating set. The lower bound is usually found by exploiting structural properties of the particular family of graphs, and it  often consists of a very technical and lengthy process (see, for example, {\cite{DH06}}). Therefore, finding good general lower bounds for the power domination number is an important problem. 
 
 An effort in that direction is the work by Stephen et al. \cite[Theorem 3.1]{SRRGW15} in which a  lower bound is presented and successfully applied to finding the power domination number of some graphs modeling chemical structures. However, their lower bound depends heavily  on the choice of a family of subgraphs satisfying certain properties. While in some graphs it is possible to find families of subgraphs that yield good lower bounds, in others it is not, and the bound depends on the family of subgraphs chosen rather than on the graph itself.
 
The lower bound for the power domination number of a hypercube presented in Dean et al. \cite{Deanetal11} is based on the following result:

\begin{thm}\label{thm:Deanetal}{\rm\cite[Lemma 2]{Deanetal11}}
	 Let $G$ be a graph with no isolated vertices, and let $S=\{u_1,\dots,u_t\}$ be a power dominating set for $G$.  Then $\Z(G)\le\sum_{i=1}^t\deg u_i.$
\end{thm}

While Dean et al. only apply Theorem \ref{thm:Deanetal} to the study of a particular family of graphs, Theorem \ref{thm:Deanetal} is the starting point in our process to obtain a general lower bound.  The next theorem, which follows from Theorem \ref{thm:Deanetal}, can be used to map zero forcing results to power dominating results and vice versa.  
 
\begin{thm}\label{Daniela}
	 Let $G$ be a graph that has an edge.  Then $\lc\frac{\Z(G)}{\Delta(G)}\rc\le \pd(G)$, and this bound is tight.
\end{thm}
\bpf   Choose a minimum power dominating set $\{u_1,\dots,u_t\}$, so $t=\pd(G)$, and observe that $\sum_{i=1}^t\deg u_i\le t\Delta(G)$.  If $G$ has no isolated vertices, the result follows from Theorem \ref{thm:Deanetal}.  Each isolated vertex of $G$ contributes one to both the zero forcing number and the power domination number, hence the result still holds.  Since $\Z(K_n)=\Delta(K_n)=n-1$ and $\pd(K_n)=1$, the bound is tight.
\epf

The next corollary is immediate from the fact that $\M(G)\le\Z(G)$.  Although weaker than Theorem \ref{Daniela}, Corollary \ref{corM} can sometimes be applied using a well known matrix such as the adjacency or Laplacian matrix of the graph, even if $\M(G)$ and $\Z(G)$ are not known. In addition, Corollary \ref{corM} permits the incorporation of a new set of tools based on linear algebra into the study of power domination.

\begin{cor}\label{corM}
	For a graph $G$ that has an edge and any matrix $A\in\SG$, $\lc\frac{\nul A}{\Delta(G)}\rc\le \pd(G)$.
\end{cor}


\subsection{Applications to computation of power domination number}\label{ssappspd}

Dorbec et al. studied the power domination problem for the tensor product of two paths \cite{DMKS}. We study the tensor product of a path and a complete graph and of a cycle and a complete graph.

\begin{prop}\label{ubPtxKn}
Let $n\geq 3$.  If $G=P_t$ with $t\geq 2$ or $G=C_t$ with $t\geq 3$, then \[\pd(G\times K_n)\leq \left\{ \begin{array}{ll}
  \lc t\over 2\rc &\mbox { if }  t\not\equiv 2 \mod 4, \\
  {t\over 2} + 1 & \mbox { if } t \equiv 2 \mod 4.
 \end{array} \right. \]
\end{prop}

\bpf
Denote the vertices of $G\times K_n$ as ordered pairs $(r,s)$ for  $1\leq r\leq t, 1\leq s \leq n$.  Define a   set $S$ in the following way (throughout, $k$ is a positive integer):

If $t=4k$, let $S = \{ (4i-2,1), (4i-1,1): 1\leq i\leq k \}$.

If $t = 4k+1$, let $S = \{ (4i-2,1), (4i-1,1): 1\leq i\leq k \} \cup \{ (4k,1) \}$.

If $t = 4k+2$, let $S = \{ (4i-2,1), (4i-1,1): 1\leq i\leq k \} \cup \{ (4k+1,1), (4k+2,1) \}$.

If $t = 4k+3$, let $S = \{ (4i-2,1), (4i-1,1): 1\leq i\leq k \} \cup \{ (4k+2,1), (4k+3,1) \} $. 

\noindent It is easy to verify that $S$ is a power dominating set for $G\times K_n$ and thus, $\pd(G\times K_n)\leq |S|$.
\epf

\begin{obs}\label{degreetensor}
The degree of an arbitrary vertex $(g,h)$ in $G\times H$ is the product $deg_G(g)deg_H(h)$. Therefore, $\Delta (G\times H)=\Delta (G)\Delta (H)$.
\end{obs}

\begin{thm}\label{PtxKn}
Let $G=P_t$ with $t\ge 2$, or $G=C_t$ with $t\geq 3$.   Suppose $t$ is odd and $n\geq t$, or suppose $t$ is even and either (1) $G=P_t$ and $n\geq {t \over 2} + 2$, or (2) $G=C_t$ and $n\geq {t \over 2}$.  Then
\[ \pd(G\times K_n)= \left\{ \begin{array}{ll}
 \lc t\over 2 \rc &\mbox { if } t\not\equiv 2 \mod 4, \\
 {t\over 2} \mbox{ or } {t\over 2} + 1 &\mbox{ if } t\equiv 2 \mod 4. 
 \end{array} \right. \]
\end{thm}

\bpf
Proposition \ref{ubPtxKn} provides an upper bound on $\pd(G\times K_n)$.  We obtain a lower bound on $\pd(G\times K_n)$ from Theorem \ref{Daniela} and results in Section \ref{szfs}, by considering two cases depending on the parity of $t$.  Observation \ref{degreetensor} yields $\Delta(G\times K_n) = \Delta(G)\Delta(K_n) = 2(n-1)$. 

If $t=2k+1$ for some positive integer $k$, then by  Theorems \ref{HCYthm}, \ref{CtKn}, and \ref{Daniela}  we know that 
\[\pd(G\times K_n)\ge \lc  {{(n-2)(2k+1)+2 }\over {2(n-1)}}\rc = \lc {{2k(n-1)+n-2k}\over {2(n-1)}} \rc = \lc  k+ {{n-2k}\over {2(n-1)}} \rc \ge k+1 \mbox{ if }n-2k>0.\] 
That is, $\lc{t\over 2}\rc\leq \pd(G\times K_n)$ if $t$ is odd and $n \ge t$.

Let $t=2k$. Define $c=0$ for $G=P_t$ and $c=2$ for $G=C_t$.  Then by  Theorems \ref{evenPtKn}, \ref{CtKn}, and \ref{Daniela} we know that 
\[\pd(G\times K_n)\ge \lc  {{(n-2)(2k)+2c}\over {2(n-1)}}\rc = \lc {{k(n-1)-k+c}\over {(n-1)}} \rc = \lc  k- {{k-c}\over {n-1}} \rc =k \mbox{ if }n-1>k-c.\] 
That is, ${t\over 2}\leq \pd(G\times K_n)$ if  $G=P_t$ and $n \ge {t \over 2}+2$, or if $G=C_t$ and $n \ge {t \over 2}$. \epf

\begin{rem}  Computations using {\em Sage} power domination software \cite{sagedata} suggest that for $t\equiv 2\mod 4$ and $n\ge 4$, the correct value is $\pd(G\times K_n)=\frac t 2 +1$.  As noted earlier, $n=3$ can behave differently, and the values computed were $\pd(P_2\times K_3)=1=\frac t 2$ and $\pd(C_6\times K_3)=3=\frac t 2$.
\end{rem}

\begin{rem}   It is shown in \cite[Theorem 4.2]{BF}\footnote{ there is a typographical error in the statement of this theorem, but it is clear from the proof that this is what is intended.} that for $m\ge n\ge 3$, \[\pd(C_n\cp C_m)\le \left\{ \begin{array}{ll}
 \lc n\over 2 \rc &\mbox { if } n\not \equiv 2 \mod 4, \\
 {n\over 2} + 1 &\mbox{ if } n\equiv 2 \mod 4. 
 \end{array} \right. \]  It follows immediately from  Theorems \ref{CnCm} and \ref{Daniela} that this inequality is an equality whenever  $n\not \equiv 2 \mod 4$, and $\frac n 2\le \pd(C_n\cp C_m)\le \frac n 2 + 1$ for $n\equiv 2\mod 4$.
There is an unpublished proof in \cite{CCC} that \[\pd(C_n\cp C_m)= \left\{ \begin{array}{ll}
 \lc n\over 2 \rc &\mbox { if } n\not \equiv 2 \mod 4, \\
 {n\over 2} + 1 &\mbox{ if } n\equiv 2 \mod 4.
 \end{array} \right.\] 
\end{rem} 

\subsection{Applications to computation of zero forcing number}\label{ssappsZ}
In the preceding section, we obtained  the power domination numbers of certain graphs 
from the corresponding zero forcing numbers.   We take the opposite approach in this section, using Theorem \ref{Daniela} and known power domination numbers to obtain the corresponding zero forcing numbers.

In \cite[Theorem 4.1]{DMKS} it was proved that:

\beq\pd(G*H)=\left\{ \begin{array}{ll}
 \gamma (G) &\mbox { if } \pd(H)=1,\\
  \gamma _{t} (G) & \mbox{ otherwise,}
 \end{array} \right.\label{lexprodeq}\eeq

\noindent where  $\gamma _t (G)$ denotes the {\em total domination number of $G$}, defined as the minimum cardinality of a dominating set $S$ in $G$ such that $N(S)=V(G)$.

Now, from Theorem \ref{Daniela} we know $\Z(G*H)\leq \pd(G*H)\Delta (G*H)$. It follows easily from the definition of lexicographic product that
$ \deg_{G*H}(g,h) = (\deg_G g)|V(H)| + \deg_H h $
for any vertex $(g,h) \in V(G*H)$, and therefore $\Delta (G*H)=\Delta (G)|V(H)|+\Delta (H)$. Then from \eqref{lexprodeq} above, we obtain

 \beq \Z(G*H)\leq \left\{ \begin{array}{ll}
 \gamma (G)\, \big(\Delta (G)|V(H)|+\Delta (H)\big) &\mbox { if } \pd(H)=1,\\
  \gamma _t (G) \, \big(\Delta (G)|V(H)|+\Delta (H)\big)& \mbox{ otherwise.}
 \end{array} \right. \label{lexprodZeq} \eeq

In particular, we obtain the following result for lexicographic products of regular graphs with low domination and power dominations numbers.

\begin{thm}\label{lex} Let $G$ and $H$ be regular graphs with degree $d_G$ and $d_H$, respectively. If $\pd(H)=1$ and $\gamma (G)=1$, then $\Z(G*H)=d_G|V(H)|+d_H$.
\end{thm}
\bpf Since $G$ is $d_G$-regular, $H$ is $d_H$-regular, and $\pd(H) = \gamma (G)=1$, the bound in \eqref{lexprodZeq} gives $\Z(G*H)\leq d_G|V(H)|+d_H$.  Moreover, since $G*H$ is $(d_G|V(H)|+d_H)$-regular, Observation \ref{deltaZ} tells us $d_G|V(H)|+d_H =\delta(G*H) \leq \Z(G*H)$.
\epf
\vspace{-10pt}

\begin{cor}\label{familieslexCnKt} For $n\ge 2$ and $m\ge 3$, $\Z(K_n*C_m)=(n-1)m + 2$.
\end{cor}


\subsection*{Acknowledgements}
This research was supported by the American Institute of Mathematics (AIM), the  Institute for Computational and Experimental Research in  Mathematics (ICERM), and the National Science Foundation (NSF) through DMS-1239280.  The authors thank AIM, ICERM, and NSF.  The authors also thank Shaun Fallat and Tracy Hall for providing information on spectra of cycles.


\end{document}